\documentclass[a4paper,12pt]{article}
\usepackage[all]{xy}
\usepackage{amsmath, amsthm}
\usepackage{amssymb,amsfonts}
\usepackage{array}
\usepackage{hyperref}

\title{Eigenvalues and Wiener index of the Zero Divisor graph $\Gamma[\mathbb {Z}_n]$}
\author{B.Surendranath Reddy,Rupali.S.Jain and N.Laxmikanth\\
\textit{surendra.phd@gmail.com},\textit{rupalisjain@gmail.com} and\\
\textit{laxmikanth.nandala@gmail.com} }{}
\date{}
\begin{document}
\theoremstyle{definition}
\newtheorem{definition}{Definition}[section]
\newtheorem{example}[definition]{Example}
\newtheorem{remark}[definition]{Remark}
\newtheorem{observation}[definition]{Observation}
\theoremstyle{plain}
\newtheorem{theorem}[definition]{Theorem}
\newtheorem{lemma}[definition]{Lemma}
\newtheorem{proposition}[definition]{Proposition}
\newtheorem{corollary}[definition]{Corollary}
\newtheorem{AMS}[definition]{AMS}
\newtheorem{keyword}[definition]{keyword}
\maketitle
\begin{abstract}
The zero divisor graph of a commutative ring $R$, denoted by $\Gamma[R]$, is a graph whose vertices are non-zero zero divisors of $R$ and two vertices are adjacent if their product is zero. In this paper, we consider the zero divisor graph $\Gamma[\mathbb{Z}_n]$ for $n=p^3$ and $n=p^2q$ with $p$ and $q$ primes. We discuss the adjacency matrix and eigenvalues of the zero divisor graph $\Gamma[\mathbb{Z}_n]$. We also calculate the Wiener index of the graph $\Gamma[\mathbb{Z}_n]$.
{\bf AMS(2010):} 13A99,05C12, 05C50.
\end{abstract}
\section{Introduction}\label{sec:1}

 The concept of the zero divisor graph of a ring $R$ was first introduced by I.Beck [1] in 1988 where he considered the set of zero divisors including zero and introduced the concepts such as diameter, grith and clique number of a zero divisor graph. Then later on Anderson and Livingston in [2], Akbari and Mohammadian in [3]  continued the study of zero divisor graph and they considered only the non-zero zero divisors. The concepts of the energy and Wiener index of a zero divisor graph was introduced by Mohammad Reza and Reza Jahani in [4],  for the ring $\Gamma[{\mathbb{Z}_n}]$ for $n=p^2$ and $n=pq$ where $p$ and $q$ are prime numbers. Motivated by the work of Mohammad Reza and Reza Jahani, in this article, we are extending the concepts of adjacency matrix,  energy and Wiener index of the zero divisor graph to the rings ${\mathbb{Z}_n}$ for $n=p^3\, \text{and}\, p^2q$, where p,q are primes.

 In this article,  section 2, is about the preliminaries and notations related to zero divisor graph of a commutative ring $R$, in section 3, we derive the standard form of the adjacency matrix of the zero divisor graph  $\Gamma[{\mathbb{Z}_n}]$, in section 4 we discuss  the eigenvalues of the graph $\Gamma[{\mathbb{Z}_n}]$ and  in section 5, we find the Wiener index of $\Gamma[{\mathbb{Z}_n}]$ for $n=p^3\, \text{and}\, p^2q$.
 \section{Preliminaries and Notations}\label{sec:2}
\begin{definition}({\textbf{Zero divisor Graph}}){\cite{1}}\\
Let R be a commutative ring with unity and $ Z[R]$ be the set of its zero divisors. Then the zero divisor graph of R denoted by  $\Gamma[R]$, is the graph(undirected) with vertex set  $Z^*[R]= Z[R]-\{\mathbf{0}\}$, the non-zero zero divisors of $R$, such that two vertices $v,w \in Z^*[R]$ are adjacent if $vw=0$.
\end{definition}
\begin{definition}({\textbf{Adjacency matrix of $\Gamma[R]$}})\\
 The adjacency matrix of the zero divisor graph $\Gamma[R]$ is the matrix $[v_{ij}]$ with rows and columns labeled by the vertices and is given by
\begin{align*} v_{ij}=\left\{
   \begin{array}{ll}
     1, & \hbox{$v_i\, \text{is adjacent to}\, v_j$;} \\
     0, & \hbox{$otherwise$.}
   \end{array} \right.
\end{align*}
Adjacency matrix of $\Gamma[R]$ is denoted by $M(\Gamma[R])$. Clearly for an undirected graph, the adjacency matrix is symmetric.
\end{definition}
\begin{definition}({\textbf{Energy of $\Gamma[R]$}}){\cite{6}}\\
The energy of the zero divisor graph $\Gamma[R]$ is defined as the sum of the absolute values of all the eigenvalues of its adjacency matrix $M(\Gamma[R])$ i.e., if $ \lambda_1, \lambda_2,\cdots, \lambda_n$ are $n$ eigenvalues of $M(\Gamma[R])$, then the energy of $\Gamma[R]$ is $E(\Gamma[R])=  \sum\limits_{i=1}^n{|\lambda_i|} $.
\end{definition}
\begin{definition}({\textbf{Weiner index of $\Gamma[R]$}}){\cite{7}}\\
 Let  $\Gamma[R]$ be a zero divisor graph with vertex set $V$. We denote the length of the shortest path between every pair of vertices $x,y\in V$ with d(x,y). Then the Wiener index of  $\Gamma[R]$ is the sum of the distances between all pair of vertices of  $\Gamma[R]$, i.e.,  $ W(\Gamma[R])= \sum\limits_{x,y\in V} d(x,y) $.
\end{definition}
\section{Adjacency Matrix of $\Gamma[\mathbb{Z}_n]$}\label{sec:3}
In this section, we derive the standard form of the adjacency matrix of $\Gamma[\mathbb{Z}_n]$ for $n=p^3$ and $n=p^2q$, $p$ and $q$ being primes. To start with we consider $p=3$ and $n=3^3$.  Then the set of non-zero zero divisors of $Z_{27}$ is $\{3,6,9,12,15,18,21,24\}$. Now by considering these non-zero zero divisors in ascending order, we get the adjacency matrix as\\ $M(\Gamma[{\mathbb{Z}_{27}}])=\begin{bmatrix}
                                                      0 & 0 & 1 & 0 & 0 & 1 & 0 & 0 \\
                                                       0 & 0 & 1 & 0 & 0 & 1 & 0 & 0 \\
                                                       1 & 1 & 1 & 1 & 1 & 1 & 1 & 1 \\
                                                       0 & 0 & 1 & 0 & 0 & 1 & 0 & 0 \\
                                                       0 & 0 & 1 & 0 & 0 & 1 & 0 & 0 \\
                                                       1 & 1 & 1 & 1 & 1 & 1 & 1 & 1 \\
                                                       0 & 0 & 1 & 0 & 0 & 1 & 0 & 0 \\
                                                       0 & 0 & 1 & 0 & 0 & 1 & 0 & 0 \\
                                                   \end{bmatrix} $\\
Instead of this, if we rearrange the non-zero zero divisors such that all the   multiples of 3 appear first and then the multiples of $3^2$, we get the adjacency matrix as\\
 $M(\Gamma[{\mathbb{Z}_{27}}])=\left[
                                                     \begin{array}{cccccc|cc}
                                                       0 & 0 & 0 & 0 & 0 & 0 & 1 & 1 \\
                                                       0 & 0 & 0 & 0 & 0 & 0 & 1 & 1 \\
                                                       0 & 0 & 0 & 0 & 0 & 0 & 1 & 1 \\
                                                       0 & 0 & 0 & 0 & 0 & 0 & 1 & 1 \\
                                                       0 & 0 & 0 & 0 & 0 & 0 & 1 & 1 \\
                                                       0 & 0 & 0 & 0 & 0 & 0 & 1 & 1 \\\hline
                                                       1 & 1 & 1 & 1 & 1 & 1 & 1 & 1 \\
                                                       1 & 1 & 1 & 1 & 1 & 1 & 1 & 1 \\
                                                     \end{array}\right]$=$\left[
\begin{array}{c|c}
O&R\\\hline
R^t&S
\end{array}
\right]$\\
where $O$ is the zero matrix of order $6$ and $R$ and $S$ are matrices of ones of size $6\times 2$ and $2$ respectively.

After getting the above interesting form of the adjacency matrix, we tried  for several values of $n=8,125,343$, etc and got the adjacency matrix in similar form and this leads us to derive the standard form for $n=p^3$ in the following theorem.
\begin{theorem}\label{r1}
Let $n=p^3$ with $p$ prime. Then the adjacency matrix of the graph $\Gamma[{\mathbb{Z}_{n}}]$ is \\
$M(\Gamma[{\mathbb{Z}_{n}}])=
 \begin{bmatrix}
                             O_{(p^2-p)} & R_{(p^2-p)\times(p-1)} \\
                                                                               R^{t}_{(p-1)\times(p^2-p)} & S_{(p-1)} \\
                          \end{bmatrix}$\\
                          where $O$ is the zero matrix  and $R$ and $S$ are the matrices of ones.
 \begin{proof}
Let $n=p^3$. Then the set of non-zero zero divisors of $\mathbb{Z}_n$ is\\
$Z^*[\mathbb{Z_\mathbf{n}}]=\{p,2p,3p,\cdots,(p^{2}-1)p\}$ with cardinality  $p^2-1$. \\
As explained in the example, we rewrite $Z^*[\mathbb{Z_\mathbf{n}}]= A\cup B $, where
\begin{align*}
 A &=\{kp|k= 1,2,\cdots,p^2-1\, \text{and}\, p\nmid k\}  \\
 B &=\{lp^2|l =1,2,3,\cdots,p-1\}
\end{align*}
and the adjacency matrix is labeled with elements of $A$ first then the elements of $B$.\\
It is clear that if $x,y \in A$ then $xy\neq 0$.\\
Since $|A|=p^2-p$, we get a zero matrix $O$ of size $p^2-p$ in the adjacency matrix.\\
And if $ x\in A, y \in B\, \text{or}\,  y \in A, x\in B$ then $xy=0$.\\
Since $|A|=p^2-p$ and $|B|=p-1$, we get a  matrix $R$ of ones  of size $(p^2-p)\times(p-1)$.\\
Also if $x,y \in B$, then $xy=0$ and hence we get a matrix $S$ of ones  of size $p-1$ in the adjacency matrix.\\
Hence the adjacency matrix by considering the elements of $A$ first and then the elements of $B$, we get\\
 $M(\Gamma\mathbb[{Z_n}])=\begin{bmatrix}
                             O_{(p^2-p)} & R_{(p^2-p)\times(p-1)} \\
                                                                               R^t_{(p-1)\times(p^2-p)} & S_{(p-1)} \\
                          \end{bmatrix}$.

 \end{proof}
\end{theorem}

Now we consider the case for $n=p^2q$ with $p$ and $q$ different primes and get the standard form of adjacency matrix in the following theorem.
\begin{theorem}\label{r2}
  Let $ n=p^2q$, where p,q are primes. Then the adjacency matrix of the graph $\Gamma[\mathbb{Z}_{n}]$ is \\
   $M(\Gamma[\mathbb{Z}_{n}]) = \left[ \begin{array}{cc|cc}
                                     O& O & R_{(p-1)(q-1)\times(p-1)} & O\\
                                     O & O& O & S_{p(p-1)\times(q-1)} \\[1mm]\hline
                                      R^{t} & O & T_{p-1}& U_{(p-1)\times (q-1)} \\
                                      O & S^{t}& U^t & O
                                  \end{array}\right]$\\
where $R,S,T$ and $U$  are matrices of ones.
 \begin{proof}
For $n=p^2q$, the non-zero zero divisors are multiples of $p$, multiples of $q$, multiples of $p^2$ and multiples of $pq$. \\
Hence $ Z^*[\mathbb{Z_\mathbf{n}}]= A\cup B\cup C\cup D $, where
\begin{align*}
A &=\{k_{1}p\,|\,k_1= 1,2,\cdots,pq-1\, \text{and}\, k_1\nmid p,q \}  \\
B &=\{k_{2}q\,|\,k_{2}=1,2,\cdots,p^2-1\, \text{and}\, k_2\nmid p\}\\
C &=\{k_{3}pq\,|\,k_{3}=1,2,\cdots,p-1\, \}\\
D &=\{k_{4}p^2\,|\,k_{4} =1,2,\cdots,q-1 \}
\end{align*}
Then $|A|=(p-1)(q-1)$, $|B|=p(p-1)$, $|C|=p-1$ and $|D|= q-1 $.\\
Since no element of $A$ is adjacent with an element of $A$ and $B$, we get the zero matrix of order $(p-1)(p+q-1)$.\\
Similarly we get zero matrices corresponding to $A$ and $D$; $B$ and $C$; and $D$ and $D$.\\
And as every element of $A$ is adjacent with every element of $C$ implies we get a matrix $R$ of ones of order $(p-1)(q-1)\times(p-1)$.\\
Similarly we get the matrix $S$ of ones corresponding to $B$ and $D$; the matrix $T$  of ones corresponding to $C$ and $C$; and the matrix $U$  of ones corresponding to $C$ and $D$.\\
Hence we get the adjacency matrix by considering the elements of $A$ first, then $B$,  then $C$ and then $D$ as\\
$M(\Gamma[\mathbb{Z}_{n}]) = \left[ \begin{array}{cc|cc}
                                     O& O & R_{(p-1)(q-1)\times(p-1)} & O\\
                                     O & O& O & S_{p(p-1)\times(q-1)} \\[1mm]\hline
                                      R^{t} & O & T_{p-1}& U_{(p-1)\times (q-1)} \\
                                      O & S^{t}& U^t & O
                                  \end{array}\right]$
 \end{proof}
\end{theorem}
\section{Eigenvalues of the zero divisor graph  $\Gamma[{\mathbb{Z}_n}]$}\label{sec:4}
In this section we discuss the non zero eigenvalues of the zero divisor graph $\Gamma[{\mathbb{Z}_n}]$ for $n=p^3$ and $n=p^2q$.
\begin{theorem}\label{r3}
Let $n=p^3$. Then the only non-zero eigenvalues of the zero divisor graph $\Gamma(\mathbb{Z}_{n})$ are given by ${\frac{(p-1)({1}\pm\sqrt{1+4p})}{2}}$.
\begin{proof} By theorem \ref{r1}, the adjacency matrix of $\Gamma(\mathbb{Z}_{p^3}))$ is given by\\ $M(\Gamma[{\mathbb{Z}_{p^3}}])=
 \begin{bmatrix}
                             O_{(p^2-p)} & R_{(p^2-p)\times(p-1)} \\
                                                                               R^{t}_{(p-1)\times(p^2-p)} & S_{(p-1)} \\
                          \end{bmatrix}$\\
                          where $O$ is the zero matrix  and $R$ and $S$ are the matrices of ones.\\
 Let $\lambda$ be a non-zero eigenvalue of $M(\Gamma[{\mathbb{Z}_{p^3}}])$.\\
 Then $|M-\lambda{I}| =\left|
 \begin{array}{cccc|cccc}
  -\lambda& 0 & \cdots & 0&1&1&\cdots&1 \\
  0 & -\lambda & \cdots & 0&1&1&\cdots&1 \\
  \vdots  & \vdots  & \ddots & \vdots &\vdots &\vdots & \vdots& \vdots\\
  0 & 0 & \cdots & -\lambda&1&1&\cdots&1\\\hline
  1&1&\cdots&1&1-\lambda&1&\cdots&1\\
  1&1&\cdots&1&1&1-\lambda&\cdots&1\\
  \vdots & \vdots& \vdots&\vdots  & \vdots & \vdots  & \ddots & \vdots\\
  1&1&\cdots&1&1&1&\cdots &1-\lambda
 \end{array}\right|= 0$\\
 Let $|M-\lambda{I}| = \begin{vmatrix}
                  A & B \\
                  C & D
                \end{vmatrix}$\\
 Then $|M-\lambda{I}|=|A||D-CA^{-1}B|=0$.\\
 Since  $A$ is a scalar matrix of order $p^2-p$, we get  $|A|=(-\lambda)^{p^2-p}$.\\
And  $ |CA^{-1}B|=
 \begin{vmatrix}
  \frac{p-p^2}{\lambda} & \frac{p-p^2}{\lambda} & \cdots & \frac{p-p^2}{\lambda} \\
  \frac{p-p^2}{\lambda} & \frac{p-p^2}{\lambda} & \cdots & \frac{p-p^2}{\lambda} \\
  \vdots  & \vdots  & \ddots & \vdots  \\
  \frac{p-p^2}{\lambda} & \frac{p-p^2}{\lambda}& \cdots & \frac{p-p^2}{\lambda}
 \end{vmatrix}$\\[2mm]
 So $ |D-CA^{-1}B|=
 \begin{vmatrix}
  (1-\lambda)+\frac{p-p^2}{\lambda} & 1+\frac{p-p^2}{\lambda} & \cdots &1+ \frac{p-p^2}{\lambda} \\
  1+\frac{p-p^2}{\lambda} &(1-\lambda)+ \frac{p-p^2}{\lambda} & \cdots &1+ \frac{p-p^2}{\lambda} \\
  \vdots  & \vdots  & \ddots & \vdots  \\
 1+ \frac{p-p^2}{\lambda} & 1+\frac{p-p^2}{\lambda}& \cdots & (1-\lambda)+\frac{p-p^2}{\lambda}
 \end{vmatrix}$\\
 Using the properties of the determinant, we get
   \begin{equation*}
  \mid D-CA^{-1}B\mid = (-1)^{p-1}\times\lambda^{p-3}[\lambda^2-(p-1)\lambda-(p^2-p)(p-1)]
 \end{equation*}
 Therefore, $|M-\lambda{I}|={\lambda^{p^2-p}}(-1)^{p-1}\times\lambda^{p-3}[\lambda^2-(p-1)\lambda-(p^2-p)(p-1)]=0$\\
 Since $\lambda\neq0$, we have $\lambda^2-(p-1)\lambda-(p^2-p)(p-1)=0 $
 \begin{align*}
\text{Hence},\,  \lambda &= { \frac{(p-1)\pm\sqrt{{(p-1}^2)+4(p^2-p)(p-1)}}{2}}\\
 &= { \frac{(p-1)\pm{(p-1)}\sqrt{{1}+4p}}{2}}\\
  &= { \frac{(p-1)({1\pm\sqrt{{1}+4p}})}{2}}
  \end{align*}
 \end{proof}
\end{theorem}
\begin{theorem}
  If p is a prime number then the Energy of the zero divisor graph $\Gamma(\mathbb{Z}_{p^3})$ is $ E(\Gamma(\mathbb{Z}_{p^3}))= (p-1)(\sqrt{{1+4p}})$.
\begin{proof}
Let $n=p^3$. Then by theorem \ref{r3}, the non-zero eigenvalues of  $\Gamma(\mathbb{Z}_{n}))$ are given by
\begin{align*}
\lambda_{1} &= { \frac{(p-1)({1+\sqrt{{1}+4p}})}{2}} ,\\ \lambda_{2} &= { \frac{(p-1)({1-\sqrt{{1}+4p}})}{2}}\\
\text {Thus the Energy}, E(\Gamma(\mathbb{Z}_{p^3}))&= |\lambda_{1}|+|\lambda_{2}|\\
&=\Big|{ \frac{(p-1)({1+\sqrt{{1}+4p}})}{2}}\Big|+\Big|{ \frac{(p-1)({1-\sqrt{{1}+4p}})}{2}}\Big|\\
&= \Big(\frac{p-1}{2}\Big)\Big(\sqrt{1+4p}+1+\sqrt{1+4p}-1\Big)\\
&=(p-1)(\sqrt{1+4p})
 \end{align*}
\end{proof}
\end{theorem}
\begin{theorem}
Let $n=p^2q$ with $p$ and $q$ primes. If $\lambda\neq0$ is a nonzero eigenvalue of the  zero divisor graph $\Gamma(\mathbb{Z}_{n})$, then  \\ $\lambda^4-(p-1)\lambda^3-2p(q-1)^{2}(p-1)\lambda^2+p(p-1)^{2}(q-1)\lambda+p(p-1)^{3}(q-1)^{2}=0 $.
\begin{proof}
By theorem \ref{r2}, the adjacency matrix is given by\\
$M(\Gamma[\mathbb{Z}_{n}]) = \left[ \begin{array}{cc|cc}
                                     O& O & R_{(p-1)(q-1)\times(p-1)} & O\\
                                     O & O& O & S_{p(p-1)\times(q-1)} \\[1mm]\hline
                                      R^{t} & O & T_{p-1}& U_{(p-1)\times (q-1)} \\
                                      O & S^{t}& U^t & O
                                  \end{array}\right]$\\
where $R,S,T$ and $U$  are matrices of ones.\\
Let $\lambda$ be a nonzero eigenvalue of $M(\Gamma\mathbb[{Z_{p^{2}q}}])$.\\ Then
$|M-\lambda{I}| =\left|
 \begin{array}{cc}
 A&B\\
 C&D
 \end{array}
 \right|= 0$\\
 where
$A=\left[ \begin{array}{cccc}
-\lambda&0&\cdots&0\\
0&-\lambda&\cdots&0\\
\vdots&\vdots&\ddots&\vdots\\
0&0&\cdots&-\lambda
\end{array}\right];\quad$
$B=C^t=\left[ \begin{array}{cccccc}
1&\cdots&1&0&\cdots&0\\
\vdots&&\vdots&\vdots&&\vdots\\
1&\cdots&1&0&\cdots&0\\
0&\cdots&0&1&\cdots&1\\
\vdots&&\vdots&\vdots&&\vdots\\
0&\cdots&0&1&\cdots&1\\
\end{array}\right]$\\
and $D=\left[ \begin{array}{cccccccc}
1-\lambda&1&\cdots&1&1&&\cdots&1\\
1&1-\lambda&\cdots&1&1&&\cdots&1\\
\vdots&\vdots&\ddots&\vdots&\vdots&&&\vdots\\
1&1&\cdots&1-\lambda&1&&\cdots&1\\
1&&\cdots&1&-\lambda&0&\cdots&0\\
1&&\cdots&1&0&-\lambda&\cdots&0\\
\vdots&&&\vdots&\vdots&\vdots&\ddots&\vdots\\
1&&\cdots&1&0&0&\cdots&-\lambda
\end{array}\right]$\\
Clearly $|A|= \lambda^ {pq+p^2-p-5}$\\[1mm]
$ |CA^{-1}B|=
 \begin{vmatrix}
  \frac{(1-p)(q-1)}{\lambda} & \cdots&  \frac{(1-p)(q-1)}{\lambda}  & 0 &\cdots & 0 \\
  \vdots&&\vdots&\vdots&&\vdots\\
  \frac{(1-p)(q-1)}{\lambda} & \cdots&  \frac{(1-p)(q-1)}{\lambda}  & 0 &\cdots & 0 \\
     0 & \cdots&0&   \frac{p(1-p)}{\lambda} & \cdots &\frac{p(1-p)}{\lambda}\\
     \vdots&&\vdots&\vdots&&\vdots\\
 0 & \cdots&0&   \frac{p(1-p)}{\lambda} & \cdots &\frac{p(1-p)}{\lambda}
  \end{vmatrix}$\\
 and hence $ |D-CA^{-1}B|=\\[1mm]
 \begin{vmatrix}
  1-\lambda+\frac{(1-p)(q-1)}{\lambda} & \cdots&  \frac{(1-p)(q-1)}{\lambda}  & 1 &\cdots & 1 \\
  \vdots&&\vdots&\vdots&&\vdots\\
  \frac{(1-p)(q-1)}{\lambda} & \cdots&  1-\lambda+\frac{(1-p)(q-1)}{\lambda}  & 1 &\cdots & 1 \\
     1 & \cdots&1&   -\lambda+\frac{p(1-p)}{\lambda} & \cdots &\frac{p(1-p)}{\lambda}\\
     \vdots&&\vdots&\vdots&&\vdots\\
 1 & \cdots&1&   \frac{p(1-p)}{\lambda} & \cdots &-\lambda+\frac{p(1-p)}{\lambda}
  \end{vmatrix}$\\
  \begin{equation*}
   = \lambda^{4}-(p-1)\lambda^{3}-2p(p-1)(q-1)\lambda^2+p(p-1)^{2}(q-1)\lambda+p(p-1)^{3}(q-1)^{2}.\\
 \end{equation*}
 Since $\lambda\neq0$, $|A|= \lambda^ {pq+p^2-p-5}$ and $|A||D-CA^{-1}B|=0$, we get\\
$ \lambda^{4}-(p-1)\lambda^{3}-2p(p-1)(q-1)\lambda^2+p(p-1)^{2}(q-1)\lambda+p(p-1)^{3}(q-1)^{2}=0$.\\
\end{proof}
\end{theorem}
\section{Wiener index of the zero divisor graph $\Gamma(\mathbb{Z}_{n})$}\label{sec:5}
In this section we calculate the Wiener index of $\Gamma(\mathbb{Z}_{n})$ for $n=p^3$ and $n=p^2q$.
\begin{theorem}
The Weiner index of the zero divisor graph  $\Gamma(\mathbb{Z}_{p^3})$ is\\ $(\frac{p-1}{2})(2p^{3}-3p-2)$.
\begin{proof}
Let $\Gamma(\mathbb{Z}_{p^3})$ be a zero divisor graph where p is a prime number.\\
Then the vertex set of  $\Gamma(\mathbb{Z}_{p^3})$ is divided into two sets
\begin{align*}
 A&=\{k_{1}p\,|\,k_{1}= 1,2,3,...,p^{2}-1\,\text{ and } k_{1}\nmid p\}  \\
  B&=\{k_{2}p^2\,|\,k_{2} =1,2,3,...,p-1\}
\end{align*}
where $|A|=p^2-p = n(say)$ and $|B|=p-1=m(say)$\\
Let $ x,y\in V(\Gamma(\mathbb{Z}_{p^3}))$\\
Then either $ x,y\in A$ or $ x,y\in B$ or $ x\in A, y\in B $ \\
Suppose  $ x,y\in A $ . Since no two elements of A are adjacent and every element of A is adjacent with an element of B, we get $d(x,y)= 2$\\
Suppose  $ x,y\in B $ . Since any two elements of B are adjacent to each other,we get $d(x,y)= 1 $\\
Suppose  $ x\in A, y\in B. $ Since every element of A is adjacent with an element of B ,we get $d(x,y)= 1 $\\
Thus
\begin{align*}
 \sum\limits_{x,y \in A} d(x,y)&= \sum\limits_{i=2}^{n} d(x_{1},x_{i}) + \sum\limits_{i=3}^{n} d(x_{2},x_{i}) + ...+  d(x_{n-1},x_{n}) \\
 &=2(n-1)+2(n-2)+...+2=n(n-1)
 \end{align*}
Similarly  $\sum\limits_{x,y \in A} d(x,y)= (m-1)+(m-2)+...+1=\frac{m(m-1)}{2}$\\
And
\begin{align*}
 \sum\limits_{x \in A, y\in B } d(x,y)&= \sum\limits_{i=1}^{m} d(x_{1},y_{i}) + \sum\limits_{i=1}^{m} d(x_{2},y_{i}) + ...+   \sum\limits_{i=1}^{m} d(x_{n},y_{i}) \\
 &=m+m+...+m(n\, times)=nm
 \end{align*}
   Therefore,the Weiner index is
 \begin{align*}
   W(\Gamma(\mathbb{Z}_{p^3}))&= \sum\limits_{x \in A, y\in B } d(x,y) + \sum\limits_{x,y \in A} d(x,y) + \sum\limits_{x,y \in B} d(x,y)\\
&= nm + n(n-1)+ \Big(\frac{m(m-1)}{2}\Big)\\
&=(p^{2}-p)(p-1)+ (p^{2}-p)(p-2)+ \frac{(p-1)(p-2) }{2} \\
&= \Big(\frac{p-1}{2}\Big)(2p^{3}-3p-2)
 \end{align*}
\end{proof}
\end{theorem}
\begin{theorem}
 The Wiener index of the zero divisor graph  $\Gamma(\mathbb{Z}_{p^2q})$ is\\ $(\frac{1}{2})[p^{2}(2p^{2}-4p-1)+p(4p^{2}q+2pq^{2}-8pq-2q+5)+2]$ \\
\begin{proof}
Let $\Gamma(\mathbb{Z}_{p^2q})$ be a zero divisor graph where p ,q are primes. \\
Then the vertex set of  $\Gamma(\mathbb{Z}_{p^2q})$ is divided into four sets
\begin{align*}
A &=\{k_{1}p\,|\,k_1= 1,2,3,...,pq-1\, \text{and}\, k_1\nmid p,q \}  \\
B &=\{k_{2}q\,|\,k_{2}=1,2,3,...,p^2-1\, \text{and}\, k_2\nmid p\}\\
C &=\{k_{3}pq\,|\,k_{3}=1,2,3,...,p-1\, \}\\
D &=\{k_{4}p^2\,|\,k_{4} =1,2,3,...,q-1 \}
\end{align*}
Then $\alpha = |A|= (p-1)(q-1)$,\\
$\beta = |B|=  p(p-1)$ ,\\
$\gamma = |C|=  (q-1)$ and\\
$\delta = |D|=  (p-1)$\\
For $ x,y\in V(\Gamma(\mathbb{Z}_{p^2q}))$, we have the following cases:
\begin{enumerate}
\item Suppose $ x,y\in A$. Since no two elements of $A$ are adjacent and every element of $A$ is adjacent with an element of $D$, we get $d(x,y)= 2$
\item If $ x,y\in B$, since no two elements of $B$ are adjacent and every element of $B$ is adjacent with an element of $C$,we get $d(x,y)= 2 $
\item If $ x,y\in C$, as no two elements of $C$ are adjacent and every element of $C$ is adjacent with an element of either $B$ or $D$ ,we get $d(x,y)= 2$
\item If $ x,y\in D $, being any two elements of D are adjacent, we get $d(x,y)= 1$
\item Now suppose $ x\in A, y\in B$, as no element of $A$ is adjacent with an element of $B$, we get $d(x,y)= 2 $\\
With similar arguments we get,
\item If $ x\in A, y\in C$,  $d(x,y)= 2 $
\item If $ x\in A, y\in D$, $d(x,y)= 1 $
\item If $ x\in B, y\in C$,  $d(x,y)= 1 $
\item If $ x\in B, y\in D$,  $d(x,y)= 2 $
\item If $ x\in C, y\in D$,  $d(x,y)= 1 $
\end{enumerate}
Hence the Weiner index is
 \begin{align*}
   W(\Gamma(\mathbb{Z}_{p^2q})&= \sum\limits_{x,y \in \{\text{same\, sets}\}} d(x,y) + \sum\limits_{x,y \in \{\text{different\, sets}\}} d(x,y)\\
&=\sum\limits_{i=2} d(x_{1},x_{i}) + \sum\limits_{i=3} d(x_{2},x_{i}) + ...+  d(x_{n-1},x_{n})\\
& +\sum\limits_{i=1} d(x_{1},y_{i}) + \sum\limits_{i=1} d(x_{2},y_{i}) + ...+ \sum\limits_{i=1} d(x_{n},y_{i})\\
 &=2[(\alpha-1)+(\alpha-2)+...+2+1] + 2[(\beta-1)+(\beta-2)+...+2+1] \\
 &+ 2[(\gamma-1)+(\gamma-2)+...+2+1]+1[(\delta-1)+(\delta-2)+...+2+1]\\
 &+2[\alpha\beta]+2[\alpha\gamma]+[\alpha\delta]+[\beta\gamma]+2[\beta\delta]+[\gamma\delta].\\
 &= \alpha(\alpha-1) + \beta(\beta-1) +\gamma(\gamma-1)+ \frac{\delta(\delta-1)}{2} \\
 &+2\alpha\beta+2\alpha\gamma+\alpha\delta+\beta\gamma+2\beta\delta+\gamma\delta
  \end{align*}
By substituting the values of $\alpha,\beta,\gamma,\delta$ in above and on solving, we get \\
$W(\Gamma(\mathbb{Z}_{p^2q}) = (\frac{1}{2})[p^{2}(2p^{2}-4p-1)+p(4p^{2}q+2pq^{2}-8pq-2q+5)+2]$.
  \end{proof}
\end{theorem}

\end{document}